\documentclass{conm-p-l}  
  
\newtheorem{theorem}{Theorem}[section]  
\newtheorem{lemma}[theorem]{Lemma}  
\newtheorem{corollary}[theorem]{Corollary}

\theoremstyle{definition}

\theoremstyle{remark}

\numberwithin{equation}{section}

\def\bC{\mathbf C}  
\def\bQ{\mathbf Q}  
  
\def\hV{\hat V}  
\def\hP{\hat P}  
  
\def\boZ{\mathbf Z}  
  
\def\bg{\mathbf g}  
  
\def\hW{\widehat W}  
  
\def\chtheta{\check\theta}  
\def\bgh{\mathbf{\hat g}}

\def\cF{{\mathcal F}}  
  
\def\cH{{\mathcal H}}

\def\cU{{\mathcal U}}

\def\ox{\otimes}  
  
\def\l{\lambda}  
\def\la{\langle}  
\def\ra{\rangle}  
\def\rkg{\ell}  
  
\begin{document}  
  
\title{A New Perspective on the Frenkel-Zhu Fusion Rule Theorem}  
  
\author{Alex J. Feingold}  
\address{Dept. of Math. Sci., The State University of New York,  
Binghamton, New York 13902-6000} 
\email{alex@math.binghamton.edu}  
 
\thanks{AJF wishes to thank the Albert Einstein Institute for 
its wonderful hospitality and support which helped complete this work.} 
  
\author{Stefan Fredenhagen}  
\address{Max-Planck-Institut f{\"u}r Gravitationsphysik,  
Albert-Einstein-Institut, M{\"u}hlenberg 1, D-14476 Golm, Germany}  
\email{Stefan.Fredenhagen@aei.mpg.de}  
  
  
\date{}  
  
  
\begin{abstract}  
In this paper we prove a formula for fusion coefficients of   
affine Kac-Moody algebras first conjectured by Walton \cite{Wal2},   
and rediscovered in \cite{Fe}. It is a reformulation of the Frenkel-Zhu affine   
fusion rule theorem \cite{FZ}, written so that it can be seen as a beautiful   
generalization of the classical Parasarathy-Ranga Rao-Varadarajan tensor product   
theorem \cite{PRV}.  
\end{abstract}  
 
\maketitle 
 
\tableofcontents  
 
\section{Introduction}  
 
Fusion rules play a very important role in conformal field theory  
\cite{Fu}, in the representation theory of vertex operator algebras  
\cite{FLM, FHL, FZ}, and in quite a few other areas. For example,  
fusion rules were used in \cite{FS} to obtain information on D-brane  
charge groups in string theory which on the other hand correspond to  
certain twisted K-groups. This line of research found a mathematical  
culmination in the theorem by Freed, Hopkins and Teleman \cite{FHT},   
showing that twisted equivariant K-theory can be identified with a  
fusion ring. In \cite{AFW} a connection was found between the fusion  
rules for the Virasoro minimal models and elementary abelian  
$2$-groups.  Further work in \cite{FW} extended this idea to find a  
connection between the fusion rules for type $A_1$ and $A_2$ on all  
levels, and elementary abelian $2$-groups and $3$-groups. This was  
extended as far as was possible in \cite{Sal1,Sal2} to the case of  
$A_\rkg$ for any rank $\rkg$ and any level.
  
In \cite{Fe} an introduction was given to the subject with major focus  
on the algorithmic aspects of computing fusion rules for affine  
Kac-Moody algebras. In particular, it was emphasized that the  
Kac-Walton algorithm \cite{Kac, Wal} for fusion coefficients is  
closely related to the Racah-Speiser algorithm for tensor product  
decompositions, which was the subject of earlier work  
\cite{F1,F2}. \cite{Fe} included a conjecture on fusion coefficients  
which restates the Frenkel-Zhu theorem \cite{FZ} in a form which shows  
it to be a beautful generalization of the classical Parasarathy-Ranga  
Rao-Varadarajan tensor product theorem \cite{PRV}. That conjecture had  
already been made by Walton \cite{Wal2} in 1994, but we believe that  
it has not been proven up until now.  
  
An outline of the organization of the paper is as follows. We give the  
definition of a fusion algebra in section two, then we give notation  
and background about finite dimensional simple Lie algebras in section  
three. This includes facts about irreducible representations,  
contravariant Hermitian forms on them, special results for $sl_2$ and  
its representations, and projection operators. In section four we  
briefly give notations about affine algebras leading to the level $k$  
fusion algebra associated with simple Lie algebra $\bg$. In section  
five we discuss tensor products of irreducible finite dimensional  
modules for $\bg$ and the PRV theorem. In section six we state the  
Frenkel-Zhu fusion rule theorem, the Walton conjecture, what it says  
in the special case when $\bg = sl_2$, and a corollary relating fusion  
coefficients to tensor product multiplicities. We begin the proof of  
the Walton conjecture by rewriting the Frenkel-Zhu theorem in several  
ways. In section seven we review the proof of the PRV theorem and  
refine it to help find a relationship between the spaces which occur  
in the Frenkel-Zhu theorem and the Walton conjecture. In section eight  
we put all these pieces together to finish the proof of the Walton  
conjecture.  
  
\section{Definition of Fusion Algebra}  
  
Let us begin with the definition of fusion algebra given by J. Fuchs  
\cite{Fu}. A fusion algebra $F$ is a finite dimensional commutative  
associative algebra over $\bQ$ with some basis  
$$B = \{x_a\ |\ a\in A\}$$   
so that the structure constants $N_{a,b}^c$ defined by  
$$x_a \cdot x_b = \sum_{c\in A} N_{a,b}^c x_c$$  
are non-negative integers. There must be a distinguished  
index $\Omega\in A$ with the following properties.  
It is required that the matrix  
$$C = [C_{a,b}] = [N_{a,b}^\Omega]$$ satisfies $C^2 = I$. Because  
$0\leq N_{a,b}^c \in\boZ$, either $C = I$ or $C$ must be an order 2  
permutation matrix, that is, there is a permutation $\sigma:A\to A$  
with $\sigma^2 = 1$ and  
$$C_{a,b} = \delta_{a,\sigma(b)}.$$  
Write $\sigma(a) = a^*$ and call $x_{a^*}$ the conjugate of $x_a$.   
Use it to define the non-negative integers  
$$N_{a,b,c} = N_{a,b}^{c^*}$$  
which, by commutativity and associativity of the algebra product, are  
completely symmetric in $a$, $b$ and $c$. Using this we also find that   
$x_\Omega$ is a multiplicative identity element in $F$ and $\Omega^* = \Omega$.  
  
In this paper we are interested in the structure constants of fusion  
algebras that are associated to affine Lie algebras.

\section{Background and Notation for Finite Dimensional Lie Algebras}  
  
Now we will introduce notations and review some basic results needed later.   
Let $\bg$ be a finite dimensional simple Lie algebra of rank $\rkg$ with  
Cartan matrix $A = [a_{ij}]$ and Cartan subalgebra $H$.   
The simple roots and the fundamental weights of $\bg$ are linear functionals  
$$\alpha_1,\cdots,\alpha_{\rkg} \quad\hbox{and}\quad\lambda_1,\cdots,\lambda_{\rkg},$$  
respectively, in the dual space $H^*$. Let the integral weight lattice $P$ be   
the $\boZ$-span of the fundamental weights, and let  
$$P^+ = \{n_1\lambda_1+\cdots+n_{\rkg}\lambda_{\rkg}\ |\  
0\leq n_1,\cdots,n_{\rkg}\in\boZ\}$$   
be the set of dominant integral weights of $\bg$, and let  
$$\theta = \sum_{i=1}^{\rkg} \theta_i \alpha_i$$  
be the highest root of $\bg$. The symmetric bilinear form $(\cdot,\cdot)$ on $H^*$  
is determined by   
$$a_{ij} = \langle\alpha_i,\alpha_j\rangle   
= \frac{2(\alpha_i,\alpha_j)}{(\alpha_j,\alpha_j)},\quad 1\leq i,j\leq \rkg$$  
and the normalization $(\theta,\theta) = 2$. The fundamental weights are   
determined by the conditions $\langle\lambda_i,\alpha_j\rangle = \delta_{ij}$  
for $1\leq i,j\leq \rkg$, and the special ``Weyl vector''  
$$\rho = \sum_{i=1}^{\rkg} \lambda_i$$  
will play an important role in several formulas. It is useful to define  
$${\check\lambda} = \frac{2\lambda}{(\lambda,\lambda)}\quad\hbox{for any }  
0\neq \lambda\in H^*,$$  
so we can write $(\lambda_i,{\check\alpha}_j) = \delta_{ij}$ and   
$a_{ij} = (\alpha_i,{\check\alpha}_j)$. We may also express  
$$\theta = \sum_{i=1}^{\rkg} \chtheta_i {\check\alpha}_i\qquad\hbox{so}\qquad  
\chtheta_i = \frac{\theta_i (\alpha_i,\alpha_i)}{2}.$$  
The Weyl group $W$ of $\bg$ is defined to be the group of endomorphisms of $H^*$  
generated by the simple reflections corresponding to the simple roots,   
$$r_i(\lambda) = \lambda - (\lambda,{\check\alpha}_i) \alpha_i,\qquad 1\leq i\leq \rkg.$$  
This is a finite group of isometries which preserve $P$.   
There is a partial order defined on $H^*$ defined by  
$$\lambda \leq \mu \quad\hbox{ if and only if }\quad \mu - \lambda =   
\sum_{i=1}^{\rkg} k_i \alpha_i \quad\hbox{for some } 0\leq k_i\in\boZ.$$  
  
For $\lambda\in P^+$ let $V^\lambda$ denote the finite dimensional irreducible  
$\bg$-module with highest weight $\lambda$. It has the weight space decomposition  
$V^\lambda = \bigoplus_{\beta\in H^*} V^\lambda_\beta$, where   
$$V^\lambda_\beta = \{v\in V^\lambda\ |\ h\cdot v = \beta(h) v, \forall h\in H\}$$  
is the $\beta$ weight space of $V^\lambda$. Of course, there are only finitely  
many $\beta\in H^*$ such that $V^\lambda_\beta$ is nonzero, and we denote by   
$\Pi^\lambda$ that finite set of such $\beta$. Since $\dim(V^\l_\l) = 1$, a  
nonzero highest weight vector $v_\l^\l\in V^\l_\l$ is determined up to a scalar.   
The dual space $(V^\lambda)^* =   
Hom(V^\lambda,\bC)$ is also an irreducible highest weight $\bg$-module, called  
the contragredient module of $V^\lambda$. The action of $\bg$ on $(V^\lambda)^*$  
is given by   
$$(x\cdot f)(v) = -f(x\cdot v) \quad\hbox{ for }\quad x\in\bg,\ f\in (V^\lambda)^*,  
\ v\in V^\lambda.$$  
The highest weight of $(V^\lambda)^*$ is denoted by $\lambda^*$, and   
equals the negative of the lowest weight of $V^\lambda$. For example, in the case   
when $\bg$ is of type $A_{\rkg}$, if $\lambda = \sum_{i=1}^{\rkg} n_i \lambda_i$   
then $\lambda^* = \sum_{i=1}^{\rkg} n_{\rkg+1-i} \lambda_i$.   
  
On $V^\lambda$ with a chosen highest weight vector, $v_\l^\l\in V^\l_\l$,    
we have a positive definite contravariant Hermitian form \cite{Kac}   
$(\cdot,\cdot): V^\l \times V^\l \to \bC$ determined by the  
following conditions: (1) $(v_\l^\l,v_\l^\l) = 1$,     
(2) For any $v,v'\in V^\l$, and any $x\in\bg$, we have  
$(x\cdot v,v') = - (v,x^\dagger\cdot v')$, where the map $x \to x^\dagger$ is the  
Chevalley involutive automorphism of $\bg$ determined by its action on the generators  
$$e_i^\dagger = -f_i, \qquad f_i^\dagger = -e_i, \qquad h_i^\dagger = -h_i,\qquad  
1\leq i\leq \rkg.$$  
Note that for any $v\in V^\l_\beta$, $v'\in V^\l_{\beta'}$, we have  
$$\beta(h_i) (v,v') = (h_i\cdot v,v') = - (v,-h_i\cdot v') = \beta'(h_i) (v,v')$$  
so $0 = (\beta - \beta')(h_i) (v,v')$ for any Cartan generator $h_i$. This means  
that if $\beta \neq \beta'$ then $(v,v') = 0$ so different weight spaces are   
orthogonal. Let $Proj_\beta^\l: V^\l \to V^\l_\beta$ denote the orthogonal projection  
operator.  
  
If $V^\l$ and $V^\mu$ are two irreducible highest weight modules with chosen highest  
weight vectors and positive definite contravariant Hermitian forms as above, then we have  
a positive definite contravariant Hermitian form on the tensor product $V^\l \ox V^\mu$   
given by  
$$(v^\l_1 \ox v^\mu_1, v^\l_2 \ox v^\mu_2) = (v^\l_1,v^\l_2)(v^\mu_1,v^\mu_2).$$  
If $V^\nu$ is an irreducible submodule of $V^\l \ox V^\mu$ then its orthogonal   
complement $(V^\nu)^\perp = \{v\in V^\l \ox V^\mu\ |\ (v,V^\nu) = 0\}$ is clearly a  
$\bg$-submodule since   
$$(x\cdot v, V^\nu) = -(v, x^\dagger\cdot V^\nu) = 0,   
\hbox{ for all }x\in\bg, v\in (V^\nu)^\perp.$$  
This shows that when the tensor product $V^\l \ox V^\mu$ is decomposed into a direct   
sum of irreducible $\bg$-modules, the distinct modules obtained are mutually orthogonal  
with respect to the contravariant Hermitian form. Let   
$Proj^{\l,\mu}_{V^\nu}: V^\l\ox V^\mu \to V^\nu$ denote the orthogonal projection  
operator from the tensor product to a particular irreducible submodule $V^\nu$.  
  
We will use certain facts about the representation theory of the simple Lie algebra   
$\bg = sl_2$ of type $A_1$ whose standard basis $\{e, f, h\}$ has the brackets   
$[h,e] = 2e$, $[h,f] = -2f$ and $[e,f] = h$.   
An irreducible finite dimensional $sl_2$-module $V^\l$ is determined by its highest  
weight, the non-negative integer $\l(h) = m$, so we write $V^\l = V(m)$.   
If $v_0$ is a highest weight vector then a basis of $V(m)$ can be written as  
$\{v_i\ |\ 0\leq i\leq m\}$ where $v_i = \frac{1}{i!} f^i v_0$ and the action of   
$\bg$ is given by the formulas:  
$$h\cdot v_i = (m - 2i) v_i, \quad f\cdot v_i = (i+1) v_{i+1}, \quad  
e\cdot v_i = (m - i + 1) v_{i-1}$$  
for $0\leq i\leq m$ with the understanding that $v_j = 0$ for $j$ outside that range.   
For any integer $p\geq 0$, we understand the operators $e^p$ and $f^p$ on $V(m)$   
to mean $p$ repetitions of the operators $e$ and $f$, respectively.   
It is easy to see that the contravariant form has values   
$(v_i,v_j) = \delta_{i,j}{m\choose i}$,  
for $0\leq i\neq j\leq m$, so the form is positive definite.  
  
\medskip  
\begin{lemma}\label{sl2lemma1}   
Let $\bg = sl_2$ and $V(m)$ be the irreducible finite dimensional $sl_2$-module with  
highest integral weight $m\geq 0$. Then for any integer $p\geq 1$, with respect to the  
positive definite contravariant Hermitian form on $V(m)$, we have an orthogonal direct   
sum decomposition  
$$V(m) = ker(f^p) \oplus Im(e^p).$$  
\end{lemma}  
  
\begin{proof} From the explicit formulas for the action, it is clear that $ker(f^p)$  
is the subspace of the $p$ lowest weight spaces with basis $\{v_{m-p+1},\cdots,v_m\}$  
and that $Im(e^p) = (ker(f^p))^\perp$ is the subspace of all other weight spaces  
with basis $\{v_0,\cdots,v_{m-p}\}$.   
\end{proof}  
\medskip  
  
We now go back to the general case of any finite dimensional simple $\bg$.   
Let $V^\l$ be an irreducible $\bg$-module, $\alpha$ any root of $\bg$, and let  
$\bg_\alpha$ be the corresponding subalgebra of $\bg$ isomorphic to $sl_2$ with   
standard basis $\{e_\alpha, f_\alpha, h_\alpha\}$.   
The Chevalley involution acts on $\bg_\alpha$ by $e_\alpha^\dagger = -f_\alpha$,   
$f_\alpha^\dagger = -e_\alpha$ and $h_\alpha^\dagger = -h_\alpha$. The complete  
reducibility of finite dimensional $sl_2$-modules gives a direct sum decomposition   
$$V^\l = \bigoplus_{i} V^\l_{\gamma_i}(m_i)$$  
into irreducible $\bg_\alpha$-modules, where $V^\l_{\gamma_i}(m_i)$ has $\bg$-highest   
weight $\gamma_i$, and $\bg_\alpha$-highest weight $\gamma_i(h_\alpha) = m_i$.   
If $V^\l_{\gamma_1}(m_1)$ is one of these, then its orthogonal complement is clearly   
a $\bg_\alpha$-submodule by the same argument as given above for the decomposition of   
a tensor product. It means that this decomposition is an orthogonal direct sum decomposition  
with respect to the contravariant Hermitian form on $V^\l$.  
  
\medskip  
\begin{lemma}\label{eflemma}   
Let $V^\l$ be an irreducible $\bg$-module, $\alpha$ any root of $\bg$, and $\bg_\alpha$   
be the corresponding subalgebra of $\bg$ isomorphic to $sl_2$. Let $\beta\in\Pi^\l$ be any  
weight of $V^\l$. Then, for any integer $p\geq 0$ such that $p+\la\beta,\alpha\ra\geq 0$,   
we have   
$$\{v\in V^\l_\beta\ |\ e_\alpha^p(v) = 0\} =   
\{v\in V^\l_\beta\ |\ f_\alpha^{p+\la\beta,\alpha\ra}(v) = 0\}.$$  
\end{lemma}  
  
\begin{proof}   
The Weyl group reflection $r_\alpha$ acts on the weights $\Pi^\l$ and  
$r_\alpha(\beta) = \beta - \la\beta,\alpha\ra \alpha$. It is also  
well-known that the operator  
$$R_\alpha = (exp(f_\alpha)) (exp((-e_\alpha)) (exp(f_\alpha)) \in  
GL(V^\l)$$   
satisfies $R_\alpha(V^\l_\beta) = V^\l_{r_\alpha(\beta)}$  
for any weight $\beta\in\Pi^\l$.  It is clear from the definition of  
$R_\alpha$ that it acts on each of the $\bg_\alpha$ submodules in the  
decomposition of $V^\l$ given in the paragraph above the lemma. For  
any $0\neq v\in V^\l_\beta$ we have $0\neq R_\alpha(v)\in  
V^\l_{r_\alpha(\beta)}$. We can write $v = \sum_i v_i$ where $v_i \in  
V^\l_{\gamma_i}(m_i)$, and $e_\alpha^p(v) = 0$ iff $e_\alpha^p(v_i) =  
0$ for each $i$, so we may assume $v$ is in one such irreducible  
$\bg_\alpha$-module. The condition $e_\alpha^p(v) = 0$ means $v$ is in one  
of the top $p$ weight spaces of its irreducible  
$\bg_\alpha$-module. This is equivalent to saying that $R_\alpha(v)$  
is in one of the bottom $p$ weight spaces, that is,  
$f_\alpha^p(R_\alpha(v)) = 0$.  
  
If $\la\beta,\alpha\ra \geq 0$ then $R_\alpha(v) = c f_\alpha^{\la\beta,\alpha\ra}(v)$ for some nonzero scalar $c$, which means   
$0 = f_\alpha^p(c f_\alpha^{\la\beta,\alpha\ra}(v)) = c f_\alpha^{p+\la\beta,\alpha\ra}(v)$.  
  
If $\la\beta,\alpha\ra < 0$ but $p+\la\beta,\alpha\ra\geq 0$ then   
$R_\alpha(v) = c e_\alpha^{-\la\beta,\alpha\ra}(v)$ for some nonzero scalar $c$ which means   
$0 = f_\alpha^p(c e_\alpha^{-\la\beta,\alpha\ra}(v)) = d f_\alpha^{p+\la\beta,\alpha\ra}(v)$  
for a nonzero scalar $d$.  
\end{proof}  
\medskip

If $V$ is any finite dimensional vector space with a positive definite Hermitian   
form and $W$ is any subspace of $V$ then $W$ has an orthogonal complement  
$W^\perp = \{v\in V \ |\ (v,w) = 0, \forall w\in W\}$   
such that $V = W \oplus W^\perp$. Let $P_W: V \to W$ be the orthogonal projection  
of $V$ onto $W$ defined by $P_W(v) = w$ where $v = w + w'$ is the unique expression  
for $v\in V$ with $w\in W$ and $w'\in W^\perp$.    
If $L:V \to V$ is any linear transformation, there is a unique adjoint  
linear transformation $L^\dagger:V \to V$ determined by the conditions   
$$(L(v),v') = (v,L^\dagger(v')), \qquad\hbox{ for all }v,v'\in V.$$  
We call $L$ self-adjoint when $L = L^\dagger$. Note that any orthogonal projection  
map is self-adjoint because if $v_1 = w_1 + w_1'$ and $v_2 = w_2 + w_2'$ for   
$w_1,w_2\in W$ and $w_1',w_2'\in W^\perp$, then   
$$(P_W(v_1),v_2) = (w_1,w_2 + w_2') = (w_1,w_2) = (w_1 + w_1',w_2) = (v_1,P_W(v_2))$$  
so $P_W^\dagger = P_W$. Also, it is clear that $P_W^2 = P_W$.   
  
Finally, later we will need the following lemma.  
\medskip  
\begin{lemma}\label{projectionlemma} Let $V = U_1 \oplus U_2$ be an orthogonal  
direct sum decomposition of a finite dimensional vector space $V$ with a positive definite  
Hermitian form, and let $W$ be any subspace of $V$. Then we have the   
orthogonal direct sum decomposition of $W$:  
$$W = P_W(U_1) \oplus (W \cap U_2).$$   
\end{lemma}  
  
\begin{proof} Let $v\in W$ be in the orthogonal complement of $P_W(U_1)$. This means  
that for any $u_1\in U_1$, we have  
$$0 = (P_W(u_1), v) = (u_1,P_W^\dagger(v)) = (u_1,P_W(v)) = (u_1,v)$$  
which means $v\in W\cap U_1^\perp = W\cap U_2$.  
\end{proof}  
  
\section{Notation for Affine Lie Algebras}  
  
Let  
$$\bgh = \bg\otimes \bC[t,t^{-1}] \oplus \bC c \oplus \bC d$$   
be the affine algebra constructed from $\bg$ with derivation $d = -t\frac{d}{dt}$   
adjoined as usual, and with Cartan subalgebra  
$$\cH = H \oplus \bC c \oplus \bC d.$$   
The simple roots and the fundamental weights of $\bgh$ are linear functionals    
$$\alpha_0, \alpha_1,\cdots,\alpha_{\rkg} \quad\hbox{and}\quad  
\Lambda_0,\Lambda_1,\cdots,\Lambda_{\rkg},$$  
respectively, in the dual space $\cH^*$. The simple roots of $\bg$ form a basis  
of $H^*$ (as do the fundamental weights), and we identify them with linear   
functionals in $\cH^*$ having the same values on $H\subseteq\cH$ and being   
zero on $c$ and $d$. Let $c^*$ and $d^*$  
in $\cH^*$ be the functionals which are zero on $H$ and which satisfy  
$$c^*(c) = 1, \quad c^*(d) = 0, \quad d^*(c) = 0, \quad d^*(d) = 1.$$  
Extend the bilinear form $(\cdot,\cdot)$ to $\cH^*$ by letting  
$$(c^*,H^*) = 0 = (d^*,H^*),\qquad (c^*,c^*) = 0 = (d^*,d^*),\ \hbox{and}\   
(c^*,d^*) = 1.$$  
Then $\alpha_0 = d^* - \theta$ and   
$$\Lambda_0 = c^*, \quad \Lambda_i   
= \theta_i \frac{(\alpha_i,\alpha_i)}{2}\ c^* + \lambda_i = \chtheta_i\ c^* + \lambda_i, 
\quad 1\leq i\leq \rkg,$$  
are determined by the conditions $\langle\Lambda_i,\alpha_j\rangle = \delta_{ij}$  
for $0\leq i,j\leq \rkg$. Let the integral weight lattice $\hP$ be   
the $\boZ$-span of the fundamental weights, and let 
$$\hP^+ = \{\sum_{i=0}^{\rkg} n_i \Lambda_i \ |\ 0\leq n_i\in\boZ\}$$  
be the set of dominant integral weights of $\bgh$.   
  
The affine Weyl group $\hW$ of $\bgh$ is the group of endomorphisms of $\cH^*$  
generated by the simple reflections corresponding to the simple roots,  
$$r_i(\Lambda) = \Lambda - (\Lambda,{\check\alpha}_i) \alpha_i,\qquad   
0\leq i\leq \rkg.$$  
This is an infinite group of isometries which preserve $\hP$.   
The canonical central element, $c\in\bgh$ acts on an irreducible $\bgh$-module   
as a scalar $k$, called the level of the module.   
We will only discuss modules with highest weight $\Lambda\in \hP^+$, which are   
the ``nicest'' in that they have affine Weyl group symmetry and satisfy the   
Weyl-Kac character formula. An irreducible highest weight $\bgh$-module is   
uniquely determined by its highest weight   
$$\Lambda = \sum_{i=0}^{\rkg} n_i \Lambda_i \in \hP^+$$  
and, if we define $\theta_0 = 1 = \chtheta_0$, then   
$$k = \Lambda(c) = \sum_{i=0}^{\rkg} n_i\Lambda_i(c)  
= \sum_{i=0}^{\rkg} n_i \theta_i \frac{(\alpha_i,\alpha_i)}{2} = 
\sum_{i=0}^{\rkg} n_i \chtheta_i.$$ 
For fixed $k$ there are only finitely many $\Lambda\in \hP^+$ with   
$\Lambda(c) = k$, and we denote that finite set by $\hP_k^+$. It is easy to   
see that $\hW$ preserves the level $k$ weights   
$\{\Lambda\in\hP\ |\ \Lambda(c) = k\}$. The affine hyperplane determined by   
the condition $\Lambda(c) = k$ can be projected onto $H^*$ and the corresponding   
action of $\hW$ is such that the simple reflections $r_i$ for $1\leq i\leq \rkg$   
act as they were defined originally on $H^*$, as isometries generating   
the finite Weyl group $W$ of $\bg$. But the new affine reflection $r_0$ acts  
as $r_0(\lambda) = \lambda - (\lambda,\theta)\theta + k\theta   
= r_{\theta}(\lambda) + k\theta$, the composition of reflection $r_\theta$   
and the translation by $k\theta$, which is not an isometry on $H^*$.   
  
Irreducible $\bgh$-modules $\hV^\Lambda$ of level $k \geq 1$ are indexed by   
$\hP_k^+$, but we can also index them by certain weights of $\bg$ as follows.   
From the formulas above we can write  
$$\Lambda = \sum_{i=0}^{\rkg} n_i \Lambda_i   
= k c^* + \sum_{i=1}^{\rkg} n_i \lambda_i.$$  
So there is a bijection between $\hP_k^+$ and the set of weights   
$\lambda = \sum_{i=1}^{\rkg} n_i \lambda_i$  
such that   
$$k = n_0 + \sum_{i=1}^{\rkg} n_i \theta_i \frac{(\alpha_i,\alpha_i)}{2}  
= n_0 + \sum_{i=1}^{\rkg} n_i \chtheta_i = n_0 + \langle\lambda,\theta\rangle.$$  
Since $n_0 \geq 0$, this is equivalent to the ``level $k$ condition''  
$$\langle\lambda,\theta\rangle   
= \sum_{i=1}^{\rkg} n_i \chtheta_i \leq k.$$  
Define the set   
$$P_k^+ = \{\lambda = \sum_{i=1}^{\rkg} n_i \lambda_i\in P^+\ |\  
\langle\lambda,\theta\rangle \leq k\}$$  
and let the index set $A$ (as in the fusion algebra definition) be $P_k^+$.   
Then we see that irreducible modules on level $k$    
correspond to $\l\in P_k^+$.   
Fix level $k\geq 1$ and write the fusion algebra product (which has not been  
defined yet!)  
$$[\lambda]\cdot [\mu] = \sum_{\nu\in P_k^+} N^{(k)\nu}_{\lambda,\mu}\ [\nu].$$  
The distinguished identity element, $[0]$, corresponds to $\Lambda = k c^*$,   
and for each $[\lambda]$ there is a distinguished conjugate $[\lambda^*]$  
such that $N^{(k)0}_{\lambda,\mu} = \delta_{\mu,\lambda^*}$.   
Knowing $N^{(k)\nu}_{\lambda,\mu}$  
is equivalent to knowing the completely symmetric coefficients  
$$N^{(k)}_{\lambda,\mu,\nu} = N^{(k)\nu^*}_{\lambda,\mu}.$$  
Let $\cF(\bg,k)$ denote this fusion algebra.  
  
\section{Tensor Product Decompositions}  
  
There is a close relationship between the product in fusion algebras  
associated with an affine Kac-Moody algebra $\bgh$ and tensor product   
decompositions of irreducible $\bg$-modules.   
Let $V^\lambda$ be the irreducible finite dimensional $\bg$-submodule  
of $\hV^\Lambda$ generated by a highest weight vector. In the special  
case when $\Lambda = k\Lambda_0 = kc^*$, that finite dimensional   
$\bg$-module is $V^0$, the one dimensional trivial $\bg$-module.   
Since $\bg$ is semisimple, any finite dimensional $\bg$-module is   
completely reducible. Therefore, we can write the tensor product   
of irreducible $\bg$-modules   
$$V^\lambda\otimes V^\mu = \sum_{\nu\in P^+}   
Mult_{\lambda,\mu}^\nu V^\nu$$  
as the direct sum of irreducible $\bg$-modules, including multiplicities.   
This decomposition is independent of the level $k$ and is part of the   
basic representation theory of $\bg$. The fusion products  
$[\lambda]\cdot [\mu]$ are obtained by a subtle truncation   
of the above summation.   
  
The Racah-Speiser algorithm gives the formula   
$$Mult_{\lambda,\mu}^\nu = \sum_{w\in W} \epsilon(w)   
Mult_\lambda(w(\nu+\rho) - \mu - \rho)$$  
where $W$ is the Weyl group of $\bg$,   
$\epsilon(w) = (-1)^{length(w)}$ is the sign of $w$,   
the Weyl vector $\rho = \sum \lambda_i$ is the sum of   
the fundamental weights of $\bg$, and   
$Mult_\lambda(\beta) = dim(V^\lambda_\beta)$ is the inner multiplicity of   
the weight $\beta$ in $V^\lambda$. Recall that   
$\Pi^\lambda = \{\beta\in H^*\ |\ dim(V^\lambda_\beta) > 0\}$ denotes the set   
of all weights of $V^\lambda$. In fact, the only weights $\nu$ for which   
$Mult_{\lambda,\mu}^\nu$ may be nonzero are those of the form   
$\nu = \beta + \mu$ where $\beta\in \Pi^\lambda$.   
  
This algorithm assumes that you can already produce the weight diagram of  
any irreducible module, $V^\lambda$, so we should have discussed that first,  
but in fact the special case of the Racah-Speiser algorithm when $\mu = 0$  
gives a recursion for the inner multiplicities of $V^\lambda$. Since $V^0$  
is the trivial one-dimensional module, $V^\lambda\otimes V^0 = V^\lambda$,   
so $Mult_{\lambda,0}^\nu = \delta_{\lambda,\nu}$ and therefore  
$$0 = \sum_{w\in W} \epsilon(w) Mult_\lambda(w(\nu+\rho) - \rho)$$  
for $\nu\neq\lambda$. One knows that $Mult_\lambda(w\lambda) = 1$  and   
$Mult_\lambda(w\nu) = Mult_\lambda(\nu)$ for all  
$w\in W$, so the above formula implies that   
$$Mult_\lambda(\nu) = - \sum_{1\neq w\in W} \epsilon(w)   
Mult_\lambda(\nu+\rho - w\rho)$$  
for $\nu\neq\lambda$. Since $\rho > w\rho$ in the partial ordering on weights,  
this gives an effective recursion for $Mult_\lambda(\nu)$.   
  
In \cite{F1,F2} Feingold studied certain patterns which occur in the tensor   
product decomposition of a fixed irreducible $\bg$-module, $V^\lambda$, with all   
other modules $V^\mu$. For fixed $\lambda$, as $\mu$ varies there are only a finite   
number of different patterns of outer multiplicities which can occur, and there  
are sets of values for $\mu$ for which the pattern is constant, called   
zones of stability for tensor product decompositions.  
We have the following precise result from \cite{F2} about when a particular weight   
$\beta$ of $V^\lambda$, reaches the zone of stability.  
  
\medskip  
\begin{theorem}\label{mythesis} Let $\lambda,\mu\in P^+$ and $\beta\in\Pi^\lambda$ be  
such that $\beta + \mu\in P^+$. Let   
$$\beta - r_{\beta,j}\alpha_j,\cdots,\beta,\cdots\beta + q_{\beta,j}\alpha_j$$  
be the $\alpha_j$ weight string through $\beta$. If  
$\langle\mu,\alpha_j\rangle \geq q_{\beta,j}$ then   
$$Mult_{\lambda,\mu}^{\beta+\mu} =   
Mult_{\lambda,\mu+\lambda_j}^{\beta+\mu+\lambda_j}.$$ \end{theorem}  
  
Since $\langle\mu+\lambda_j,\alpha_j\rangle = \langle\mu,\alpha_j\rangle + 1$,  
it is clear that $\langle\mu,\alpha_j\rangle \geq q_{\beta,j}$ implies  
$$Mult_{\lambda,\mu}^{\beta+\mu} =   
Mult_{\lambda,\mu+m\lambda_j}^{\beta+\mu+m\lambda_j}\qquad\hbox{for all }m\geq 1.$$  
This result shows that for fixed $\lambda\in P^+$ and fixed  
$\beta\in\Pi^\lambda$, the tensor product multiplicities   
$Mult_{\lambda,\mu}^{\beta+\mu}$ have zones of stability as $\mu$ varies, and   
it is sufficient to study the finite number of $\mu$ such that   
$\langle\mu,\alpha_j\rangle \leq q_{\beta,j}$ for $1\leq j\leq \rkg$.  
  
There is another important result about tensor product coefficients  
which played a role in \cite{F1,F2}. In 1977 Prof.\ Bertram Kostant  
drew the attention of Feingold to the following beautiful result of  
Parthasarathy, Ranga Rao and Varadarajan \cite{PRV}, which is here  
rewritten slightly.  
  
\medskip  
\begin{theorem}\label{prv} \cite{PRV} Let $\lambda,\mu\in P^+$ and   
$\beta\in\Pi^\lambda$ be such that $\beta + \mu\in P^+$. Let $\rkg = rank(\bg)$   
and let $0\neq e_j\in\bg_{\alpha_j}$ be a root vector corresponding  
to the simple root $\alpha_j$ for $1\leq j\leq \rkg$. Then   
$$Mult_{\lambda,\mu}^{\beta+\mu} =   
dim\{v\in V^\lambda_\beta\ |\ e_j^{\langle\mu,\alpha_j\rangle+1} v = 0,   
1\leq j\leq \rkg\}.$$\end{theorem}  
  
\section{The Frenkel-Zhu Theorem and its reformulation}  
  
Now let us turn to the Frenkel-Zhu fusion rule theorem for affine Kac-Moody  
algebras. (Note that this is closely related to results of Gepner-Witten   
\cite{GW}, which appeared much earlier in the physics literature. Also,  
see Haisheng Li \cite{Li}.)  
  
\begin{theorem}\label{fz1} \cite{FZ} Let $\lambda,\mu,\nu\in P^+_k$, and  
let $0\neq e_\theta\in\bg_\theta$ be a root vector of $\bg$ in the $\theta$ root space of  
$\bg$. Let $v_\nu^\nu \in V^\nu$ be a highest weight vector and write  
$$\cH' = Hom_\bg(V^\l \ox V^\mu \ox V^\nu,\bC).$$  
Then the level $k$ fusion coefficient $N_{\l,\mu,\nu}^{(k)}$, which is completely  
symmetric in $\l$, $\mu$ and $\nu$, equals the dimension of the vector space   
$$FZ_k(\l,\mu,\nu) = \{f\in \cH'\ |\   
f(e_\theta^{k-\langle\nu,\theta\rangle+1} V^\l \ox V^\mu \ox v^\nu_\nu) = 0\}.$$  
\end{theorem}  
\medskip  
  
We now state the main result of this paper, the theorem, conjectured by Walton,   
which is a blending of the PRV and FZ theorems, showing that the FZ theorem is   
actually a beautiful generalization of the PRV theorem.   
  
\medskip  
\begin{theorem} \label{AJFThm} For $\lambda,\mu\in P^+_k$, $\beta\in\Pi^\lambda$  
such that $\beta+\mu\in P^+_k$, we have $N_{\lambda,\mu}^{(k)\ (\beta+\mu)}$  
equals the dimension of the space   
$$W_k^+(\lambda,\beta,\mu) =   
\{v\in V^\lambda_\beta\ |\ e_j^{\langle\mu,\alpha_j\rangle+1} v = 0,   
1\leq j\leq\rkg, \hbox{ and } e_\theta^{k-\langle\beta+\mu,\theta\rangle+1} v  
= 0\}.$$ \end{theorem}  
\medskip  
  
In \cite{Wal2} the statement of the conjecture is slightly different from above, with  
the condition $e_\theta^{k-\langle\beta+\mu,\theta\rangle+1} v = 0$ replaced by the condition   
$f_\theta^{k-\langle\mu,\theta\rangle+1} v = 0$. The equivalence of these two conditions  
is precisely the content of Lemma \ref{eflemma}.  
  
Theorem \ref{AJFThm} implies the following result, which tells the level $k$ at which   
the fusion coefficient associated with a single weight $\beta\in\Pi^\lambda$   
equals the tensor product multiplicity associated with that weight.  
  
\medskip  
\begin{corollary} \label{AJFCor} Suppose $\lambda,\mu\in P^+_k$, and  
$\beta\in\Pi^\lambda$ is such that $\beta+\mu\in P^+_k$. Let the $\theta$   
weight string through $\beta$ in $\Pi^\lambda$ be   
$\beta - r\theta,\cdots,\beta,\cdots,\beta + q\theta$. Then   
$k\geq \la\mu,\theta\ra + r$ implies   
$N_{\lambda,\mu}^{(k)\ (\beta+\mu)} = Mult_{\lambda,\mu}^{\beta+\mu}$.  
\end{corollary}  
\medskip  
  
Before starting the proof of the theorem, we will show how it  
reproduces the well known fusion coefficients in the special  
case when $\bg = sl_2$, where $\rkg = 1$, $\theta = \alpha_1$, and  
$P_k^+ = \{n_1\l_1\ |\ n_1\in\boZ, 0\leq n_1\leq k\}$. In this case we  
use the notation $[n_1]$ instead of $n_1\l_1$, so $V^{[n_1]} = V(n_1)$  
is the irreducible $\bg$-module with highest weight $[n_1]$. The  
weights of $V^{[n_1]}$ are $\{\beta = [n_1 - 2i]\ |\ 0\leq i\leq n_1  
\}$ and each weight space $V^{[n_1]}_{[n_1-2i]}$ is  
one-dimensional. For $0\leq n_1\leq n_2\in\boZ$, the tensor product  
decomposition  
$$V^{[n_1]} \ox V^{[n_2]} = \bigoplus_{i=0}^{n_1} V^{[n_1 + n_2 -  
2i]}$$   
is well-known. If $[n_1],[n_2]\in P_k^+$ then the fusion  
product corresponds to a truncation of this tensor product, so that  
only terms $[n_1 + n_2 - 2i]\in P_k^+$ could appear, with coefficients  
no larger than $1$. Note that the following Corollary  
\ref{sl2corollary} says the truncation is somewhat stronger than that,  
requiring $n_1 + n_2 - 2i \leq k - i$.  Since there is a symmetry  
between $n_1$ and $n_2$, it is not surprising to also find the  
condition $i\leq n_2$ symmetric to the assumption $i\leq n_1$.  
  
\medskip  
\begin{corollary} \label{sl2corollary} For $0\leq n_1,n_2\leq k$,   
$0\leq i\leq n_1$ with $0\leq n_1 - 2i + n_2 \leq k$, the $sl_2$ fusion coefficient  
$N_{[n_1],[n_2]}^{(k)\ [n_1+n_2-2i]}$ equals $1$ if   
$i \leq n_2$ and $n_1 + n_2 - 2i \leq k - i$, zero otherwise.  
\end{corollary}  
  
\begin{proof}For $1\leq i\leq n_1$, the raising   
operator $e_1 = e_\theta$ sends $V^{[n_1]}_{[n_1-2i]}$ isomorphically onto   
$V^{[n_1]}_{[n_1-2i+2]}$, and kills the highest weight space $V^{[n_1]}_{[n_1]}$.   
This means that for $v\in V^{[n_1]}_{[n_1-2i]}$ and $p\geq 0$,   
$$e_1^{p+1} v = 0\quad\hbox{ iff } \quad n_1 < n_1 - 2i + 2(p+1)    
\quad\hbox{ iff } \quad i \leq p.$$  
The conditions on $v$ in the Walton space   
$$W_k^+([n_1],[n_1-2i],[n_2]) =   
\{v\in V^{[n_1]}_{[n_1-2i]}\ |\ e_1^{n_2+1} v = 0   
 \hbox{ and } e_1^{k-(n_1+n_2-2i)+1} v = 0\}$$  
are then $i\leq n_2$ and $n_1 + n_2 - 2i \leq k - i$. When these are satisfied, we have  
$$W_k^+([n_1],[n_1-2i],[n_2]) = V^{[n_1]}_{[n_1-2i]}$$   
so the $sl_2$ fusion coefficient $N_{[n_1],[n_2]}^{(k)\ [n_1+n_2-2i]} = 1$, and otherwise,   
it is zero.  
\end{proof}  
\medskip  
  
In order to prove Theorem \ref{AJFThm} we must understand the connection between   
the PRV theorem, the statement of the theorem and the FZ theorem. We begin by  
rewriting the FZ theorem in a slightly different form.   
We can define a $\bg$-module map   
$$\Phi: Hom(V^\lambda \ox V^\mu, V^{\nu^*}) \to   
Hom(V^\lambda \ox V^\mu \ox V^\nu,\bC)$$  
by   
$$(\Phi f)(v^\lambda \ox v^\mu \ox v^\nu) = (f(v^\lambda \ox v^\mu))(v^\nu).$$  
It is easy to check that this is a $\bg$-module map and an isomorphism. In   
general, for $V$ and $W$ any two $\bg$-modules, $Hom(V,W)$ is a $\bg$-module  
under the action, $(x\cdot L)(v) = x\cdot(L(v)) - L(x\cdot v)$ for any $v\in V$  
and any $L\in Hom(V,W)$. It may be helpful to use the notations  
$\pi_V: \bg\to End(V)$, $\pi_W: \bg\to End(W)$, and $\pi: \bg\to End(Hom(V,W))$  
to distinguish the representations of $\bg$ on these three spaces. Then the   
above equation is saying that $\pi(x)(L) = \pi_W(x)\circ L - L\circ\pi_V(x)$.   
  
We also have the definition of the space of $\bg$-module maps from $V$ to $W$,   
\begin{align}  
Hom_\bg(V,W) &= \{L\in Hom(V,W)\ |\ \pi(x)(L) = 0, \forall x\in\bg\} \notag \\  
&= \{L\in Hom(V,W)\ |\ \pi_W(x)\circ L = L\circ\pi_V(x), \forall x\in\bg\}. \notag  
\end{align}  
If $v\in V_\beta$ is a weight vector of weight $\beta$, that is, for any $h\in H$,   
$\pi_V(h)v = \beta(h)v$, and $L$ is any $\bg$-module map, then $\pi_W(h)L(v) =   
L(\pi_V(h)v) = L(\beta(h)v) = \beta(h)L(v)$ shows that $L(V_\beta)\subseteq W_\beta$.   
If $Proj_\beta^V:V\to V_\beta$ and $Proj_\beta^W:W\to W_\beta$ are the orthogonal  
projection operators, then it is easy to see that $L(Proj_\beta^V(v)) = Proj_\beta^W(L(v))$  
for any $v\in V$.  
  
Since $\Phi$ is a $\bg$-module isomorphism, it is clear that it restricts to  
an isomorphism  
$$\Phi: Hom_\bg(V^\lambda \ox V^\mu, V^{\nu^*}) \to  
Hom_\bg(V^\lambda \ox V^\mu \ox V^\nu,\bC).$$  
We wish to describe the preimage of the space $FZ_k(\l,\mu,\nu)$ under   
$\Phi$. Since $\Phi$ is an isomorphism, $f\in FZ_k(\l,\mu,\nu)$ is of the form  
$\Phi g$ for a unique element $g\in Hom_\bg(V^\l \ox V^\mu, V^{\nu^*})$.   
The conditions on $f$ mean that   
$$(g(e_\theta^{k-\langle\nu,\theta\rangle+1} V^\l \ox V^\mu))  
(v^\nu_\nu) = 0.$$  
This allows us to rewrite the FZ theorem as follows.  
  
\medskip  
\begin{theorem}\label{fz2} \cite{FZ} Let $\l,\mu,\nu\in P^+_k$, and  
let $0\neq e_\theta\in\bg_\theta$ be a root vector of $\bg$ in the   
$\theta$ root space of $\bg$. Let $v_\nu^\nu \in V^\nu$ be a highest weight vector and write  
$$\cH = Hom_\bg(V^\l \ox V^\mu,V^{\nu^*}).$$  
Then the level $k$ fusion coefficient $N_{\l,\mu,\nu}^{(k)}$ equals the dimension of   
the space  
\begin{equation}  
FZ'_k(\l,\mu,\nu) = \{g\in \cH\ |\   
g(e_\theta^{k-\la\nu,\theta\ra+1} V^\l \ox V^\mu)(v_\nu^\nu) = 0\}.   
\label{fz2:eq1}  
\end{equation}  
\end{theorem}  
\medskip  
  
There is a natural isomorphism of $\bg$-modules   
\begin{equation}\label{defofpsi}  
\Psi : Hom(V^*,W) \to W \ox V  
\end{equation}  
which is defined as follows. For any $L\in Hom(V^*,W)$,  
$$\Psi(L) = \sum_{j=1}^d L(v_j^*) \ox v_j$$  
where $d = dim(V) = dim(V^*)$, $\{v_1,\cdots,v_d\}$ is any basis of $V$ and  
$\{v_1^*,\cdots,v_d^*\}$ is the dual basis of $V^*$, that is, the basis such  
that $v_i^*(v_j) = \delta_{ij}$. The inverse map sends a basic tensor   
$w\ox v \in W \ox V$ to the element in $Hom(V^*,W)$ which sends any  
$f\in V^*$ to $f(v)w\in W$. We will always choose the basis of $V$ to consist  
of weight vectors, and if $v_j$ has weight $\mu_j$, so that for any $h\in H$,  
$\pi_V(h) v_j = \mu_j(h) v_j$, then it is easy to see that the weight of  
the dual vector $v_j^*$ is $-\mu_j$. Namely, by the definition of the representation of   
$\bg$ on the dual space $V^*$, for $1\leq i\leq d$ we have  
\begin{align}  
(\pi_{V^*}(h) v_j^*)(v_i) &= - v_j^*(\pi_V(h) v_i) = - v_j^*(\mu_i(h) v_i) =  
- \mu_i(h) v_j^*(v_i) \notag\\  
&= - \mu_i(h) \delta_{ij} = - \mu_j(h) \delta_{ij} = - \mu_j(h) v_j^*(v_i)  
\notag\end{align}  
which says that $\pi_{V^*}(h) v_j^* = - \mu_j(h) v_j^*$. So $\Pi^{\lambda^*} =  
- \Pi^\lambda$. This means that a highest weight vector $v_\nu^\nu\in V_\nu^\nu$ has a dual lowest weight vector $v_{-\nu}^{\nu^*}\in V_{-\nu}^{\nu^*}$, and all other weight vectors of $V^{\nu^*}$ with weights   
above $-\nu$ are zero on $v_\nu^\nu$. In other words, with respect to the positive definite Hermitian form on the irreducible module $V^{\nu^*}$, the orthogonal complement of the lowest weight space $V_{-\nu}^{\nu^*}$ is the subspace of linear functionals in $V^{\nu^*}$ that send $v_\nu^\nu$ to $0$.  
We now see that   
\begin{equation}  
FZ'_k(\l,\mu,\nu) = \{g\in\cH\ |\   
g(e_\theta^{k-\la\nu,\theta\ra+1} V^\l \ox V^\mu)\in   
(V_{-\nu}^{\nu^*})^\perp \}.\label{fz2eq3}  
\end{equation}  
  
For any $g\in\cH$ we know that $Im(g)$ is a submodule of $V^{\nu^*}$, so if $g\neq 0$ then $g$ is surjective.   
Also, $g$ sends weight vectors to weight vectors of the same weight, and $g$ sends highest (resp., lowest)  
weight vectors to highest (resp., lowest) weight vectors. $V^{\nu^*}$ has a one dimensional highest  
weight space in which we have chosen a basis vector $v_{\nu^*}^{\nu^*}\in V_{\nu^*}^{\nu^*}$.   
$V^{\nu^*}$ also has a one dimensional lowest weight space in which we have chosen a basis vector   
$v_{-\nu}^{\nu^*}\in V_{-\nu}^{\nu^*}$. The tensor product $V^\l \ox V^\mu$ decomposes into the direct sum of irreducible modules, but $g$ must send any highest (resp., lowest) weight vector whose weight is not $\nu^*$   
(resp., not $-\nu$) to zero, so it sends all irreducible components whose highest weight is not $\nu^*$ to zero.   
The dimension of the space of highest (resp., lowest) weight vectors in $V^\l \ox V^\mu$ of weight $\nu^*$ 
(resp., $-\nu$) is the tensor product multiplicity $M = Mult_{\l,\mu}^{\nu^*}$, so we may choose a basis   
$\{u_1,\cdots,u_M\}$ of that HWV space $U^+$ (resp., LWV space $U^-$) and determine $g_i\in\cH$ uniquely by the conditions $g_i(u_j) = \delta_{i,j} v_{\nu^*}^{\nu^*}$ (resp., $g_i(u_j) = \delta_{i,j} v_{-\nu}^{\nu^*}$)   
for $1\leq i,j\leq M$. Then $\{g_1,\cdots,g_M\}$ is a basis of $\cH$. Let us denote by $\cU(\bg)$ the universal enveloping algebra of $\bg$. It is clear that $g_i$ takes the submodule   
$\cU(\bg)u_i$ isomorphically to $V^{\nu^*}$ and sends all other irreducible submodules $\cU(\bg)u_j$, $j\neq i$,   
of the tensor product to zero, so it is essentially an orthogonal projection from the tensor product to one of its components followed by an isomorphism. Let $Proj^{\l,\mu}_{U^+}$ be the orthogonal projection from $V^\l \ox V^\mu$ to the subspace  
of highest weight vectors of weight $\nu^*$, and let $Proj^{\l,\mu}_{U^-}$ be the orthogonal projection from   
$V^\l \ox V^\mu$ to the subspace of lowest weight vectors of weight $-\nu$. Then for any $v\in V^\l \ox V^\mu$,  
write $v = u + v' + v''$ where $u = Proj^{\l,\mu}_{U^-}(v)\in U^-$, $v'$ is of weight $-\nu$ but is orthogonal   
to $U^-$ so is not a lowest weight vector and must be a sum of vectors from irreducible components whose highest weights are not $\nu^*$, and $v''$ is a sum of vectors of weights not $-\nu$. Then $g(v) = g(u) + g(v') + g(v'')$   
with $g(u)\in V^{\nu^*}_{-\nu}$, and $g(v') = 0$ and $g(v'')$ is a sum of vectors of weights not $-\nu$, so   
$Proj^{\nu^*}_{-\nu}(g(v)) = g(u) = g(Proj^{\l,\mu}_{U^-}(v))$. A similar argument applies to $U^+$, so we have shown that for any $g\in\cH$ we have  
\begin{align}  
g\circ Proj^{\l,\mu}_{U^+} &= Proj^{\nu^*}_{\nu^*} \circ g\ , \\  
g\circ Proj^{\l,\mu}_{U^-} &= Proj^{\nu^*}_{-\nu} \circ g .  
\end{align}  
But this means that we can rewrite the Frenkel-Zhu space in (\ref{fz2eq3}) as   
\begin{align}  
FZ'_k(\l,\mu,\nu) &= \{g\in\cH\ |\ Proj^{\nu^*}_{-\nu} g(e_\theta^{k-\la\nu,\theta\ra+1} V^\l\ox V^\mu) = 0\} \label{fz2eq4} \\  
&= \{g\in\cH\ |\ g(Proj^{\l,\mu}_{U^-} (e_\theta^{k-\la\nu,\theta\ra+1} V^\l\ox V^\mu)) = 0\}.\label{fz2eq5}  
\end{align}  
  
\section{Review of the proof of the PRV theorem}  
  
Now we will review the proof of the PRV theorem and see if it allows us to find an isomorphism between the   
Frenkel-Zhu space $FZ'_k(\l,\mu,\nu)$ and the Walton space   
$W_k^+(\lambda,\beta,\mu)$ when $\nu^* = \beta + \mu$.  
  
In the proof of the PRV theorem one looks at the $\bg$-module $V =  
Hom(V^{\mu^*}, V^\l)$, where $\pi:\bg\to End(V)$ denotes the  
representation. As noted above (see eq.\ (\ref{defofpsi})), $V \cong  
V^\l \ox V^\mu$, and this isomorphism is given by the map $\Psi$ which  
sends irreducible components in $V$ to isomorphic irreducible  
components in $V^\l \ox V^\mu$.  The proof begins by considering the  
subspace of all lowest weight vectors (LWVs) in $V$,  
$$U = \{L\in V\ |\ \pi(f_i)L = 0,\ 1\leq i\leq \rkg\}$$  
where $\rkg = rank(\bg)$ and $e_i$, $f_i$, $h_i$ are the generators of $\bg$  
with the usual Serre relations. Then   
$$L\in U \quad\hbox{ iff }\quad \pi_\l(f_i)\circ L = L\circ \pi_{\mu^*}(f_i),\quad\hbox{for  }1\leq i\leq \rkg.$$  
It is clear that $U$ is invariant under the operators  
$\pi(h_j)$, so it has a weight space decomposition   
$$U = \bigoplus_{m=1}^r U_m$$  
where $U_m = \{L\in U\ |\ \pi(h)L = -\nu_m(h)L,\ \forall h\in H\}$ is the  
$-\nu_m$-weight space, $-\nu_1, \cdots, -\nu_r$ are the distinct lowest   
weights of irreducible components in $V$ whose corresponding highest weights   
are $\nu_1^*,\cdots,\nu_r^*$. Furthermore,   
$$dim(U_m) = Mult_{\lambda,\mu}^{\nu_m^*}$$  
is the multiplicity of $V^{\nu_m^*}$ in the tensor product $V^\lambda \ox V^\mu$  
because the independent vectors in $U_m$ each generate a distinct irreducible  
component in $V$. Let $v_1^* = v_{\mu^*}^{\mu^*}$ be a highest weight   
vector (HWV) in $V^{\mu^*}$ of weight $\mu^*$ dual to $v_1 =   
v_{-\mu^*}^{\mu}$ a LWV in $V^{\mu}$ of weight $-\mu^*$.   
The key step in the proof of  
the PRV theorem is the following lemma.  
  
\medskip  
\begin{lemma} \label{prvlemma}Define the linear map $\xi : U \to V^\l$ by   
$$\xi(L) = L(v_1^*),\qquad \forall L\in U.$$  
Then $\xi$ is injective and the range of $\xi$ equals  
$$V' = \{v\in V^\l\ |\ \pi_\l(f_i)^{\la\mu^*,\alpha_i\ra+1} v = 0,\   
1\leq i\leq \rkg\}.$$  
\end{lemma}  
  
\begin{proof} Because the highest weight vector $v_1^*\in V^{\mu^*}$ satisfies  
$$\pi_{\mu^*}(f_i)^{\la\mu^*,\alpha_i\ra+1} v_1^* = 0$$  
for $1\leq i\leq \rkg$, we have   
$$\pi_\l(f_i)^{\la\mu^*,\alpha_i\ra+1} L(v_1^*) =   
L(\pi_{\mu^*}(f_i)^{\la\mu^*,\alpha_i\ra+1} v_1^*) = 0$$  
so $\xi(U) \subseteq V'$. Let $\bg = \bg^- \oplus H\oplus \bg^+$ be the  
triangular decomposition of $\bg$, where $\bg^-$ is the Lie subalgebra of $\bg$   
generated by the negative root vectors, that is, the span of    
$f_1, \cdots, f_{\rkg}$ and all their multibrackets, and similarly $\bg^+$ is  
generated by the positive root vectors. Let $\cU(\bg)$ be the  
universal enveloping algebra of $\bg$ and extend the meaning of any representation  
of $\bg$ to include the representation of the associative algebra $\cU(\bg)$. We  
may also have use for the universal enveloping algebras $\cU(\bg^-)$ and   
$\cU(\bg^+)$. It is well known that $\cU(\bg^-)$ is spanned by all products of   
the form $y = f_{i_1} \cdots f_{i_s}$ for any $s\geq 0$ and any $1\leq i_j\leq \rkg$   
for $1\leq j\leq s$, and that $V^{\mu^*} = \cU(\bg^-) v_1^*$ is spanned by all   
vectors of the form   
$$\pi_{\mu^*}(y) v_1^*   
= \pi_{\mu^*}(f_{i_1}) \cdots \pi_{\mu^*}(f_{i_s}) v_1^*$$  
for $y$ as above. If $L(v_1^*) = 0$ for some $L\in U$ then we get  
$$0 = \pi_\l(y) L(v_1^*) = L(\pi_{\mu^*}(y) v_1^*)$$  
showing that $L = 0$ and therefore $\xi$ is injective. Let $v\in V'$ be  
arbitrary and try to define $L\in V$ by  
$$L(\pi_{\mu^*}(y) v_1^*) = \pi_\l(y) v$$  
for any $y\in \cU(\bg^-)$. If $\pi_{\mu^*}(y) v_1^* = 0$ then it is known  
that $y$ can be written  
$$y = \sum_{i=1}^\rkg y_i\ f_i^{\la\mu^*,\alpha_i\ra+1}$$  
for some $y_i \in \cU(\bg^-)$, so $\pi_\l(y) v = 0$. This means that $L$ is  
well-defined on $\pi_{\mu^*}(\cU(\bg^-)) v_1^* = V^{\mu^*}$. By the  
definition of the linear map $L$ we have  
$$(\pi_\l(f_i)\circ L)(\pi_{\mu^*}(y) v_1^*) = \pi_\l(f_i y) v =   
L(\pi_{\mu^*}(f_i y) v_1^*)   
= (L\circ \pi_{\mu^*}(f_i))(\pi_{\mu^*}(y) v_1^*)$$  
which shows that $\pi_\l(f_i)\circ L = L\circ \pi_{\mu^*}(f_i)$ so $L\in U$.   
This completes the argument that $\xi$ is an isomorphism from $U$ to $V'$.   
\end{proof}  
\medskip  
  
Now suppose that $L\in U_m$ for some $1\leq m\leq r$, so $\pi(h)L = -\nu_m(h)L$  
for any $h\in H$. But $\pi(h)L = \pi_\l(h)\circ L - L\circ\pi_{\mu^*}(h)$  
so $\xi(L)\in V^\l_{\mu^*-\nu_m}$ has weight $\mu^*-\nu_m$ because  
\begin{align}  
\pi_\l(h)(Lv_1^*) &= L(\pi_{\mu^*}(h) v_1^*) - \nu_m(h)Lv_1^* =   
L(\mu^*(h) v_1^*) - \nu_m(h)Lv_1^* \notag\\  
&= (\mu^* - \nu_m)(h) Lv_1^*.\notag \end{align}  
This shows that $\xi$ provides an isomorphism between each subspace $U_m$ and  
$$V'_{\mu^*-\nu_m} = \{v\in V^\l_{\mu^*-\nu_m}\ |\   
\pi_\l(f_i)^{\la\mu^*,\alpha_i\ra+1} v = 0,\ 1\leq i\leq \rkg\}.$$  
The PRV notation for this subspace is $V^-(\l;\mu^*-\nu_m,\mu^*)$ and  
their result is the formula for the tensor product multiplicity  
$$Mult_{\lambda,\mu}^{\nu_m^*} = dim(V^-(\l;\mu^*-\nu_m,\mu^*)).$$  
Replacing $f_i$ by $e_i$ in the definition of the space  
$V^-(\l;\gamma,\mu^*)$ one gets another space,  
$$V^+(\l;\gamma,\mu^*) = \{v\in V^\l_\gamma\ |\   
\pi_\l(e_i)^{\la\mu^*,\alpha_i\ra+1} v = 0,\ 1\leq i\leq \rkg\}.$$   
In the proof of the PRV theorem it is shown that   
$$dim(V^-(\l;\gamma,\mu^*)) = dim(V^+(\l;-\gamma^*,\mu))$$  
by using an automorphism coming from the longest element of the Weyl group, $W$.   
Then the final result of the PRV theorem is that  
$$Mult_{\l,\mu}^{\nu_m^*} = dim(V^+(\l;\nu_m^*-\mu,\mu)).$$  
  
To understand this we must discuss the longest element and a little bit of the  
theory of Lie groups. First it is necessary to know that the elements of the  
Weyl group are in one-to-one correspondence with the Weyl chambers in $H^*$.   
The dominant chamber, $P^+$, corresponding to the identity element in $W$, is  
also associated with a choice of simple roots,   
$\Delta = \{\alpha_1,\cdots,\alpha_\rkg\}$, or with a choice of positive roots,  
$R^+$, by the condition $\lambda\in P^+$ iff $\la\lambda,\alpha_i\ra \geq 0$,  
for $1\leq i\leq \rkg$. The opposite chamber $-P^+$ defined by the conditions   
$\la\lambda,\alpha_i\ra \leq 0$ is related to $P^+$ by a unique element   
$w_0\in W$ such that $w_0(P^+) = -P^+$, which means $w_0(\Delta) = -\Delta$, and   
$w_0(R^+) = R^-$. This is the longest element whose length is the number of   
positive roots and whose order is 2. For example, in type $A_2$,   
$w_0 = r_1 r_2 r_1 = r_\theta$, but in type $B_2$,   
$w_0 = r_1 r_2 r_1 r_2 \neq r_\theta$. Since $w_0(\Delta) = -\Delta$, there   
is an order 2 permutation $\sigma\in S_\rkg$ such that   
$w_0(\alpha_i) = -\alpha_{\sigma(i)}$ for $1\leq i\leq \rkg$.   
If $\nu\in P^+$ then $w_0(\nu) = -\nu^*$ is the lowest weight in $\Pi^\nu$,   
so we have  
$$\la\nu,\alpha_i\ra = \la w_0(\nu),w_0(\alpha_i)\ra =   
\la -\nu^*,-\alpha_{\sigma(i)}\ra = \la\nu^*,\alpha_{\sigma(i)}\ra.$$  
We use $\nu^* = -w_0(\nu)$ to extend the definition of dual weight to any  
$\nu\in H^*$. Note that $\theta$ is the highest weight of the adjoint   
representation and $-\theta = w_0(\theta) = -\theta^*$ is the lowest weight,   
so $\theta^* = \theta$. Therefore, for any $\nu\in H^*$ we have  
$$\la\nu,\theta\ra = \la w_0(\nu),w_0(\theta)\ra =   
\la -\nu^*,-\theta\ra = \la \nu^*,\theta\ra.$$  
  
We say $\pi_V:\bg\to End(V)$ is an integrable representation when $\pi_V(H)$  
acts diagonalizably on $V$ and all $\pi_V(e_i)$ and $\pi_V(f_i)$ are locally  
nilpotent on $V$. This is certainly true for $V$ any finite dimensional  
$\bg$-module, including the adjoint representation, $\bg$ itself, so that   
$exp(\pi_V(x))\in GL(V)$ and $exp(ad(x))\in Aut(\bg)$ for all $x = e_i$,   
$x =  f_i$ and $x = h\in H$. It is not hard to check that  
$$(exp(\pi_V(x)))\ \pi_V(y)\ (exp(\pi_V(x)))^{-1} = \pi_V(exp(ad(x))y)$$  
for all $y\in\bg$. Of particular interest are the elements  
$$r_i^{\pi_V} = (exp(\pi_V(f_i))) (exp(\pi_V(-e_i))) (exp(\pi_V(f_i))) \in   
GL(V)$$  
for $1\leq i\leq \rkg$. It is known \cite{Kac} that   
$r_i^{\pi_V}(V_\mu) = V_{r_i(\mu)}$ for any weight $\mu$ of $V$, and   
$r_i^{ad}(\bg_\alpha) = \bg_{r_i(\alpha)}$ for any root $\alpha$ of $\bg$.   
If the longest element is written as a product of simple reflections,   
$w_0 = r_{i_1} \cdots r_{i_s}$, then we have corresponding elements  
$$w_0^{\pi_V} = r_{i_1}^{\pi_V} \cdots r_{i_s}^{\pi_V} \in GL(V)\quad  
\hbox{ and }\quad  
w_0^{ad} = r_{i_1}^{ad} \cdots r_{i_s}^{ad} \in Aut(\bg)$$  
such that   
$$w_0^{\pi_V} \circ \pi_V(y) \circ (w_0^{\pi_V})^{-1} = \pi_V(w_0^{ad}(y))$$  
so using $y = h\in H$ we can get  
$$w_0^{\pi_V}(V_\mu) = V_{w_0(\mu)}\quad\hbox{ and }\quad  
w_0^{ad}(\bg_\alpha) = \bg_{w_0(\alpha)}.$$  
In particular, this means that for $1\leq i\leq \rkg$, we have  
$$w_0^{ad}(e_i)\in\bg_{w_0(\alpha_i)} = \bg_{-\alpha_{\sigma(i)}}$$  
so $w_0^{ad}(e_i) = c_i f_{\sigma(i)}$ for some $0\neq c_i\in\bC$  
and $w_0^{ad}(f_i) = c_i^{-1} e_{\sigma(i)}$. Then we have   
$$w_0^{\pi_V} \circ \pi_V(f_i) = \pi_V(w_0^{ad}(f_i))\circ w_0^{\pi_V}  
= c_i^{-1} \ \pi_V(e_{\sigma(i)})\circ w_0^{\pi_V}$$  
and for any power, $p_i$,   
$$w_0^{\pi_V} \circ \pi_V(f_i)^{p_i}   
= c_i^{-p_i}\ \pi_V(e_{\sigma(i)})^{p_i} \circ w_0^{\pi_V}.$$   
Using $p_i = \la\mu^*,\alpha_i\ra+1$ and $V = V^\lambda$, we see that   
$w_0^{\pi_V}$ provides an isomorphism between   
$$V^-(\lambda;\gamma,\mu^*) = \{v\in V^\lambda_\gamma\ |\   
\pi_\lambda(f_i)^{\la\mu^*,\alpha_i\ra+1} v = 0,\ 1\leq i\leq\rkg\}$$  
and   
$$V^+(\lambda;-\gamma^*,\mu) = \{v\in V^\lambda_{-\gamma^*}\ |\   
\pi_\lambda(e_i)^{\la\mu,\alpha_i\ra+1} v = 0,\ 1\leq i\leq\rkg\}.$$  
Since $w_0^{ad}(\bg_\theta) = \bg_{-\theta}$ we also have  
$w_0^{ad}(e_\theta) = c\ f_\theta$ for some $0\neq c\in\bC$  
and  for any power, $p$,   
$$w_0^{\pi_V} \circ \pi_V(f_\theta)^p  
= c^{-p}\ \pi_V(e_\theta)^{p} \circ w_0^{\pi_V}.$$  
Applying $w_0^{\pi_\lambda}$ to the space $W_k^+(\lambda,\beta,\mu)$ in   
Theorem \ref{AJFThm} gives the isomorphic space   
\begin{align} &W_k^-(\lambda,-\beta^*,\mu^*) = \notag\\  
&\{v\in V^\l_{-\beta^*}\ |\ \pi_\l(f_j)^{\la\mu^*,\alpha_j\ra+1} v = 0,   
1\leq j\leq \rkg, \hbox{ and } \pi_\l(f_\theta)^{k-\la\beta+\mu,\theta\ra+1} v  
= 0\}. \label{walton2}\end{align}  
It is clear that $W_k^-(\l,-\beta^*,\mu^*)$ is a subspace of  
$V^-(\l;-\beta^*,\mu^*)$,    
$$W_k^-(\l,-\beta^*,\mu^*) = \{v\in V^-(\l;-\beta^*,\mu^*)\ |\   
\pi_\l(f_\theta)^{k-\la\beta+\mu,\theta\ra+1} v = 0\}$$  
which corresponds by $\xi$ to a subpace of $U$. Our next step  
is to find the condition on $L\in U$ which corresponds to this subspace.  
  
\section{Conclusion of the proof}  
  
The root vector $f_\theta\in\bg_{-\theta}$ can be expressed as some  
multibracket of the simple root vectors $f_1,\cdots,f_\rkg$, so $L\in U$  
implies that $\pi(f_\theta)L = 0$ so $\pi_\l(f_\theta)\circ L =   
L\circ \pi_{\mu^*}(f_\theta)$. Furthermore, since $-\theta$ is the lowest  
root of $\bg$, $[f_\theta,f_i] = 0$ for $1\leq i\leq\rkg$, so in any   
representation of $\bg$, the representatives of these root vectors commute.   
For any $p\geq 1$ define the subspace of $V'$  
$$V'(p) = \{v\in V^\l\ |\ \pi_\l(f_i)^{\la\mu^*,\alpha_i\ra+1} v = 0,\   
1\leq i\leq\rkg, \pi_\l(f_\theta)^p v = 0\}.$$  
Then for any $L\in U$, $\xi(L)\in V'(p)$ iff   
$\pi_\l(f_\theta)^p L(v_{\mu^*}^{\mu^*}) = 0$ iff  
$\pi_\l(y)\pi_\l(f_\theta)^p L(v_{\mu^*}^{\mu^*}) = 0$ for all $y\in \cU(\bg^-)$.  
But since $\pi_\l(y)$ commutes with $\pi_\l(f_\theta)$, and since   
$\pi_\l(y) L(v_{\mu^*}^{\mu^*}) = L(\pi_{\mu^*}(y) v_{\mu^*}^{\mu^*})$ and  
$\cU(\bg^-) v_{\mu^*}^{\mu^*} = V^{\mu^*}$, so   
$$\xi(L)\in V'(p) \quad\hbox{iff}\quad \pi_\l(f_\theta)^p L(V^{\mu^*}) = 0  
\quad\hbox{iff}\quad L(V^{\mu^*})\subseteq Ker(\pi_\l(f_\theta)^p).$$  
Then $\xi$ provides an isomorphism from the subspace    
$$U(p) = \{L\in U\ |\ \pi_\l(f_\theta)^p L(V^{\mu^*}) = 0\} =   
\{L\in U\ |\ L(V^{\mu^*})\subseteq Ker(\pi_\l(f_\theta)^p) \} $$   
to $V'(p)$. Let $-\nu$ be one of the weights $-\nu_m$ which occur in the   
weight space decomposition of $U$, corresponding to a highest weight module  
$V^{\nu^*}$ where $\nu^* = \beta + \mu$ so   
$\la\beta + \mu,\theta\ra = \la\nu^*,\theta\ra = \la\nu,\theta\ra$.   
We have seen that $\xi$ provides an   
isomorphism between $U_{-\nu}$ and $V'_{\mu^*-\nu} = V^-(\l;\mu^*-\nu,\mu^*) =  
V^-(\l;-\beta^*,\mu^*)$, so it also provides an isomorphism between   
corresponding weight spaces  
$$U_{-\nu}(p) = \{L\in U_{-\nu}\ |\ \pi_\l(f_\theta)^p L(V^{\mu^*}) = 0\}  
= \{L\in U_{-\nu}\ |\ L(V^{\mu^*})\subseteq Ker(\pi_\l(f_\theta)^p) \}$$   
and   
$$V'_{-\beta^*}(p) = \{v\in V^\l_{-\beta^*}\ |\   
\pi_\l(f_i)^{\la\mu^*,\alpha_i\ra+1} v = 0,\ 1\leq i\leq\rkg,   
\pi_\l(f_\theta)^p v = 0\},$$  
which will equal the Walton space $W_k^-(\l,-\beta^*,\mu^*)$ when   
$p = k - \la\nu,\theta\ra+1$.   
  
\medskip  
\begin{lemma}\label{finallemma} For any integer $p\geq 1$ we have  
$$\Psi(U_{-\nu}(p)) = (Ker(\pi_\l(f_\theta)^p) \ox V^\mu) \cap \Psi(U_{-\nu})$$  
and we have the orthogonal direct sum decomposition  
$$\Psi(U_{-\nu}) = \Psi(U_{-\nu}(p)) \oplus   
Proj^{\l,\mu}_{\Psi(U_{-\nu})} (Im(\pi_\l(e_\theta)^p) \ox V^\mu).$$  
\end{lemma}  
  
\begin{proof} Apply the isomorphism $\Psi$ to $U_{-\nu}(p)$ to get the subspace   
$$\Psi(U_{-\nu}(p)) = \{\Psi(L)\in V^\l \ox V^\mu\ |\ L\in U_{-\nu}(p)\}$$  
of certain lowest weight vectors of weight $-\nu$ in $V^\l \ox V^\mu$.   
Recall the definition   
$$\Psi(L) = \sum_{j=1}^d L(v_j^*) \ox v_j$$  
where $d = dim(V^\mu) = dim(V^{\mu^*})$, $\{v_1,\cdots,v_d\}$ is a basis of $V^\mu$  
and $\{v_1^*,\cdots,v_d^*\}$ is the dual basis of $V^{\mu^*}$. Then we see that  
$$\Psi(L) \in Ker(\pi_\l(f_\theta)^p) \ox V^\mu, \hbox{ for all }L\in U_{-\nu}(p)$$  
since $L(v_j^*)\in Ker(\pi_\l(f_\theta)^p)$ for $1\leq j\leq d$. Of course,   
$\Psi(L) \in \Psi(U_{-\nu})$, so we get containment in one direction. Now suppose that   
$\Psi(L) \in \Psi(U_{-\nu})$ and $\Psi(L) \in Ker(\pi_\l(f_\theta)^p) \ox V^\mu$, so  
for $1\leq j\leq d$ we have $L(v_j^*)\in Ker(\pi_\l(f_\theta)^p)$, giving   
$L\in U_{-\nu}(p)$ so $\Psi(L) \in \Psi(U_{-\nu}(p))$.   
  
Let $\bg_\theta \cong sl_2$ be the subalgebra with basis $e_\theta$, $f_\theta$ and  
$h_\theta = [e_\theta,f_\theta]$. As mentioned in Section 3,   
$V^\l$ has a decomposition into the orthogonal direct  
sum of irreducible $\bg_\theta$-modules,   
$$V^\l = \bigoplus_{i} V^\l_{\gamma_i}(m_i)$$  
where $dim(V^\l_{\gamma_i}(m_i)) = m_i + 1$ and the highest weight of $V^\l_{\gamma_i}(m_i)$ is  
$\gamma_i\in\Pi^\l$ so $m_i = \gamma_i(h_\theta)$. 
Also recall from Section 3 that from the representation theory of   
$sl_2$, on each irreducible component we have the orthogonal decomposition  
$$V^\l_{\gamma_i}(m_i) = Ker(\pi_\l(f_\theta)^p) \oplus Im(\pi_\l(e_\theta)^p)$$  
into the $p$ lowest $h_\theta$ weight spaces and the rest. So we also get the   
orthogonal decomposition  
$$V^\l = Ker(\pi_\l(f_\theta)^p) \oplus Im(\pi_\l(e_\theta)^p).$$  
Of course, in the first equation above we mean the kernel and image of those operators   
restricted to each irreducible component. This gives an orthogonal decomposition   
$$V^\l\ox V^\mu = Ker(\pi_\l(f_\theta)^p)\ox V^\mu \oplus Im(\pi_\l(e_\theta)^p)\ox V^\mu.$$  
Lemma \ref{projectionlemma} applied to this decomposition of the tensor product gives the  
orthogonal direct sum decomposition of the subspace $\Psi(U_{-\nu})$ as stated.  
\end{proof}  
  
Let $\{\Psi(L_1),\cdots,\Psi(L_{d_p})\}$ be a basis of the first summand  
$\Psi(U_{-\nu}(p))$ in the above decomposition of $\Psi(U_{-\nu})$, and let   
$\{\Psi(L_{d_p+1}),\cdots,\Psi(L_M)\}$ be a basis of the second summand, where   
$M = Mult_{\l,\mu}^{\nu^*} = dim(U_{-\nu}) = dim(\Psi(U_{-\nu}))$. Then there is a basis,   
$\{g_1,\cdots,g_{d_p},\cdots,g_M\}$ of $\cH = Hom_\bg(V^\l\ox V^\mu,V^{\nu^*})$   
determined by the conditions $g_i(\Psi(L_j)) = \delta_{i,j} v^{\nu^*}_{-\nu}$ for   
$v^{\nu^*}_{-\nu}$ a lowest weight vector in $V^{\nu^*}$. The subspace   
$$\cH(K_p) = \{g\in \cH\ |\ g(\Psi(U_{-\nu}(p))) = 0 \}$$  
of elements of $\cH$ that vanish on the first summand,   
has basis $\{g_{d_p+1},\cdots,g_M\}$ and the subspace   
$$\cH(I_p) = \{g\in \cH\ |\ g(Proj^{\l,\mu}_{\Psi(U_{-\nu})} (Im(\pi_\l(e_\theta)^p) \ox V^\mu)) = 0 \}$$  
of elements of $\cH$ that vanish on the second summand,  
has basis $\{g_1,\cdots,g_{d_p}\}$ so $d_p = dim(\cH(I_p))$. Remember that the dimension  
of the Walton space $W_k^-(\l,-\beta^*,\mu^*)$ is $d_p$ when $p = k - \la\nu,\theta\ra + 1$.   
But in that case, $\cH(I_p)$ equals the Frenkel-Zhu space   
\begin{equation}  
FZ'_k(\l,\mu,\nu)   
= \{g\in\cH\ |\ g(Proj^{\l,\mu}_{\Psi(U_{-\nu})} (e_\theta^{k-\la\nu,\theta\ra+1} V^\l\ox V^\mu)) = 0\}\notag  
\end{equation}  
so we have completed the proof of Theorem \ref{AJFThm}.


\begin{thebibliography}{BKMW}  
  
\bibitem[AFW]{AFW}  
F. Akman, A. Feingold, M. Weiner,  
{\em Minimal model fusion rules from 2-groups,}  
Lett. Math. Phys.  {\bf  40} (1997),  159--169.  
  
\bibitem[F1]{F1} A. J. Feingold  
{\em Zones of uniform decomposition in tensor products,}  
Proc. Amer. Math. Soc. {\bf 70} (1978), 109--113.  
  
\bibitem[F2]{F2} A. J. Feingold  
{\em Tensor products of finite dimensional modules for complex semisimple  
Lie algebras,} Lie Theories and Their Applications, Proceedings of the 1977 Annual  
Seminar of the Canadian Mathematical Congress, Queen's Papers in Pure and   
Applied Mathematics, No. 48,   
(A. J. Coleman and P. Ribenboim, eds.), Queen's University, Kingston, Ontario,  
1978, pp. 394--397.  
  
\bibitem[Fe]{Fe}  
A. Feingold , {\em Fusion rules for affine Kac-Moody algebras},   
Kac-Moody Lie Algebras and Related Topics, Ramanujan International Symposium on Kac-Moody Algebras  
and Applications, Jan. 28-31, 2002, Ramanujan Institute for Advanced Study in Mathematics,   
University of Madras, Chennai, India, N. Sthanumoorthy, Kailash Misra, Editors,   
Contemporary Mathematics {\bf 343}, American Mathematical Society, Providence, RI, 2004, 53--96.  
  
\bibitem[FW]{FW}  
A. Feingold, M. Weiner,  
{\em Type A fusion rules from elementary group theory,}  
Proceedings of the Conference on Infinite-Dimensional Lie Theory   
and Conformal Field Theory, S. Berman, P. Fendley, Y. Huang, K. Misra,   
and B. Parshall, Editors, Contemporary Mathematics, Vol. 297, Amer. Math. Soc,  
Providence, RI, 2002.  
  
\bibitem[FS]{FS} S. Fredenhagen, V. Schomerus,  
{\em  Branes on group manifolds, gluon condensates, and twisted K  
theory},  
JHEP 04 (2001), 007.   
  
\bibitem[FHT]{FHT} D. S. Freed, M. J. Hopkins, C. Teleman,  
{\em Twisted K-theory and loop group representations},   
arXiv:math/0312155  

\bibitem[FHL]{FHL}  I. B. Frenkel, Yi-Zhi Huang, J. Lepowsky,  
{\em  On axiomatic approaches to vertex operator algebras and modules},  
Memoirs Amer. Math. Soc., 104, No. 594, Amer. Math. Soc., Providence, RI, 1993.   
  
\bibitem[FLM]{FLM}   
I. B. Frenkel, J. Lepowsky, A. Meurman,  
{\em  Vertex Operator Algebras and the Monster},  
Pure and Applied Math., 134, Academic Press, Boston, 1988.  
  
\bibitem[FZ]{FZ}  I. B. Frenkel, Y. Zhu,  
{\em Vertex operator algebras associated to representations of   
affine and Virasoro algebras}, Duke Math. J. {\bf  66 } (1992),  123--168.   
  
\bibitem[Fu]{Fu} J. Fuchs,  
{\em Fusion rules in conformal field theory}, Fortsch. Phys. {\bf  42 } (1994),  1--48.  
  
\bibitem[GW]{GW} D. Gepner, E. Witten  
{\em String theory on group manifolds}, Nuclear Physics {\bf B278} (1986), 493--549.  
  
\bibitem[Kac]{Kac}  V. G. Kac,   
{\em  Infinite Dimensional Lie Algebras}, Cambridge University Press,  
Third Edition, Cambridge, 1990.   
  
\bibitem[Li]{Li}  Haisheng Li,  
{\em  The regular representation, Zhu's $A(V)$-theory and induced modules},  
Journal of Algebra 238 (2001), no. 1, 159--193.  
  
\bibitem[PRV]{PRV} K. R. Parthasarathy, R. Ranga Rao, V. S. Varadarajan,  
{\em Representations of complex semi-simple Lie groups and Lie algebras},   
Annals of Mathematics, 2nd Series, Vol. 85 (1967), 383--429.  
  
\bibitem[Sal1]{Sal1} Omar Saldarriaga,  
{\em Fusion Algebras, Symmetric Polynomials, Orbits of N-Groups, and Rank-Level Duality},   
Ph.D. Dissertation, State University of New York, Binghamton, New York, 2004, 
arXiv:math.RA/0406303.  
  
\bibitem[Sal2]{Sal2} Omar Saldarriaga,  
{\em Fusion Algebras, Symmetric Polynomials, and $S_k$-orbits of $\boZ_N^k$},   
Journal of Algebra 312 (2007), 257--293.  
  
\bibitem[Wal]{Wal}  M. A. Walton,  
{\em Algorithm for WZW fusion rules: a proof},  
Phys. Lett. {\bf B241} (1990), No. 3, 365--368.  
  
\bibitem[Wal2]{Wal2} M. A. Walton,  
{\em Tensor products and fusion rules},  
Canadian Journal of Physics {\bf 72} (1994), 527--536.  
  
\end{thebibliography}
\end{document}